\documentclass[12pt]{article} 

\usepackage{amsmath,amsxtra,amssymb,latexsym, amscd}
\usepackage[mathscr]{eucal}

\begin{document}

\title{ \huge\bf Beurling spectrum of functions 
\\ in Banach space
}

\def\1{\rule{0cm}{0cm}} \def\qd{\rule{3mm}{3mm}} \def\BB{$\bullet$}
\renewcommand{\arraystretch}{1.25}
\renewcommand{\theequation}{\thesectn.\arabic{equation}}
\def\sce{\setcounter{equation}{0}}  \newcounter{sectn} 
\newcounter{sbsect}
\def\sect#1{\addtocounter{section}{1}\sce\setcounter{sbsect}{0}%
        \renewcommand{\thesectn}{\thesection}\1\smallskip\\
        {\1\hspace{-2em}\large\bf\thesectn.\qquad #1\smallskip\par}}
\def\subsect#1{\addtocounter{sbsect}{1}\sce%
        
\renewcommand{\thesectn}{\thesection:\Alph{sbsect}}\1\smallskip\\
        {\bf\1\hspace{-1.5em}\thesectn.\qquad #1\smallskip\par}}
\newtheorem{Theorem}{THEOREM} \newtheorem{Lemma}[Theorem]{LEMMA}
\newtheorem{Corollary}[Theorem]{COROLLARY}
\def\thm#1#2{\be{Theorem}{\lb{#1} #2}} \def\LEM#1#2{\BE{Lemma}{\LB{#1} 
#2}}
\def\COR#1#2{\BE{Corollary}{\LB{#1} #2}}
\def\proof{\bigskip\noindent {\sc Proof:}\qquad}
\def\REM{\1\smallskip\par\noindent{\bf REMARK:}\qquad }
\def\qed{\hfill$\quad$\qd\medskip\\} \def\ds{\displaystyle}
\def\LB#1{\label{#1}} \def\BE#1#2{\begin{#1} #2 \end{#1}}
\def\EQ#1#2{\BE{equation}{\LB{#1} #2}} \def\ARR#1#2{\BE{array}{{#1} 
#2}}
\def\DES#1{\BE{description}{#1}} \def\QT#1{\BE{quote}{#1}}
\def\ENUM#1{\BE{enumerate}{#1}} \def\ITM#1{\BE{itemize}{#1}}
 \def\COM#1{\par\noindent{\bf COMMENT:\quad\sl #1}\par\noindent}
\def\mapsfrom{\hbox{$\;{\leftarrow}\kern-.15em{\mapstochar}\:\:$}}
\def\vv{\kern.344em{\rule[.18ex]{.075em}{1.32ex}}\kern-.344em}
\def\RE{\mbox{\rm I\kern-.21em R}} \def\CX{\mbox{\rm \vv C}}
\def\imp{\Rightarrow} \def\emb{\hookrightarrow} 
\def\wk{\rightharpoonup}
\def\rd{\dot{\1}} \def\d{\cdot} \def\+{\oplus} \def\x{\times}
\def\<{\langle} \def\>{\rangle} \def\o{\circ} \def\at#1{\Bigr|_{#1}}
\def\cd{\partial} \def\grad{\nabla} \def\L{\left} \def\R{\right}
\def\bx{\mathbf{x}} \def\by{\mathbf{y}} \def\bS{\mathbf{S}}
\def\I{{\cal I}} \def\A{{\cal A}} \def\D{{\cal D}}\def\bc{\mathbf{c}}
\def\bx{\mathbf{x}} \def\by{\mathbf{y}} \def\bS{\mathbf{S}}
\def\H{{\mathcal H}} \def\U{{\cal U}} \def\D{{\cal 
D}}\def\bc{\mathbf{c}}
\def\eq{equation} \def\de{differential \eq} \def\pde{partial \de}
\def\sol{solution} \def\pb{problem} \def\bdy{boundary} 
\def\fn{function}
\def\dde{delay \de} \def\ev{eigenvalue}
\def\R{\mathbb R}
\def\C{\mathbb C}
\author{
{\bf Dang Vu Giang}\\
Hanoi Institute of Mathematics\\
18 Hoang Quoc Viet, 10307 Hanoi, Vietnam\\
{\footnotesize          e-mail: $\<$dangvugiang@yahoo.com$\>$}\\
\1\\
}
\maketitle 
{\footnotesize
\noindent {\bf Abstract.}    We are interested in Beurling spectrum of  $\mathbb X-$valued functions with application in functional delay differential equations.

\medskip


\medskip

}

\eject

\sect{\bf Beurling spectrum of $\mathbb{X}-$valued functions  and  
 differential operator}
\bigskip

\par\noindent  
In this paper, $\cal S$ denotes the set of Schwartz functions and $\phi ,\varphi ,\psi $ denote Schwartz functions. Let
\[\hat{\varphi }\left( s \right)={\cal F}\left( \varphi ,s \right)=\int\limits_{-\infty }^{\infty }{{{e}^{-ist}}\varphi \left( t \right)dt}\]
 denote the Fourier transform of $\varphi .$ Then $\hat{\varphi }={\cal F}\left( \varphi  \right)$ is also a Schwartz function and we have the inversion formula
\[\varphi \left( s \right)={{\cal F}^{-1}}\left( \hat{\varphi },s \right)=\frac{1}{2\pi }\int\limits_{-\infty }^{\infty }{{{e}^{ist}}\hat{\varphi }\left( t \right)dt}.\]
 Let $\left( \mathbb{X},{{\left\| {\cdot} \right\|}_{\mathbb{X}}} \right)$ denote a complex Banach space and let $BC\left( \mathbb{R}\to \mathbb{X} \right)$ denote the set of all $\mathbb{X}$-valued bounded continuous functions $u:\mathbb{R}\to \mathbb{X}.$ For a given function $u\in BC\left( \mathbb{R}\to \mathbb{X} \right),$ let ${{\left\| u \right\|}_{\infty }}=\sup \left\{ {{\left\| u\left( t \right) \right\|}_{\mathbb{X}}}:t\in \mathbb{R} \right\}.$ Then $\left( BC\left( \mathbb{R}\to \mathbb{X} \right),{{\left\| {\cdot} \right\|}_{\infty }} \right)$ itself is a Banach space. We define the derivative $Du=\dot{u}$ of $u\in BC\left( \mathbb{R}\to \mathbb{X} \right),$ as usual,
\[Du\left( s \right)=\dot{u}\left( s \right)=\underset{\delta \to 0}{\mathop{\lim }}\,\frac{u\left( s+\delta  \right)-u\left( s \right)}{\delta }\quad\hbox{ 
(if exists)}.\]
The differential operator is linear but unbounded on $\left( BC\left( \mathbb{R}\to \mathbb{X} \right),{{\left\| {\cdot} \right\|}_{\infty }} \right)$. But note that for every $\lambda \in \mathbb{C}\backslash i\mathbb{R}$ the operator $\lambda -D$ is invertible. More exactly,
\[{{\left( \lambda -D \right)}^{-1}}u\left( \xi  \right)=\left\{ \begin{matrix}
   \int\limits_{0}^{\infty }{{{e}^{-\lambda t}}u\left( \xi +t \right)dt\quad\text{   if  Re}\left( \lambda  \right)>0}  \\
   -\int\limits_{0}^{\infty }{{{e}^{\lambda t}}u\left( \xi -t \right)dt\quad\text{   if  Re}\left( \lambda  \right)<0.}  \\
\end{matrix} \right.\]
Hence,  the spectrum of the differential operator is $i\mathbb{R},$ and in notation, Spec$\left( D \right)=i\mathbb{R}.$ Clearly, the inverse of $\lambda -D$ is bounded operator on $\left( BC\left( \mathbb{R}\to \mathbb{X} \right),{{\left\| {\cdot} \right\|}_{\infty }} \right).$ Moreover, 
\[\int\limits_{-\infty }^{\infty }{\varphi \left( t \right){{\left( \lambda -D \right)}^{-n}}u\left( t \right)dt}=\int\limits_{-\infty }^{\infty }{\left[ {{\left( \lambda -D \right)}^{-n}}\varphi \left( t \right) \right]u\left( t \right)dt}\]
 for any $u\in $ $BC\left( \mathbb{R}\to \mathbb{X} \right).$ The convolution $\varphi *u$ of $u$ with a Schwartz function  is defined by letting
\[ \varphi *u\left( s \right)=\int\limits_{-\infty }^{\infty }{\varphi \left( s-t \right)u\left( t \right)dt.}\]
 Clearly, $\varphi *u\in BC\left( \mathbb{R}\to \mathbb{X} \right).$ The Fourier transform of the convolution of  two  functions $\varphi $ and $\psi \in {{L}^{1}}\left( \mathbb{R} \right)$  is ${\cal F}\left( \varphi *\psi  \right)=\hat{\varphi }\hat{\psi }$ and consequently, $\varphi *\psi ={{\cal F}^{-1}}\left( \hat{\varphi }\hat{\psi } \right)$ and ${{\cal F}^{-1}}\left( \varphi \psi  \right)={{\cal F}^{-1}}\left( \varphi  \right)*{{\cal F}^{-1}}\left( \psi  \right).$ Also ${\cal F} \left( \varphi \psi  \right)={{\left( 2\pi  \right)}^{-1}}{\cal F} \left( \varphi  \right)*{\cal F} \left( \psi  \right)$ and ${{\left\| \varphi *\psi  \right\|}_{\infty }}\le {{\left\| \varphi  \right\|}_{1}}{{\left\| \psi  \right\|}_{\infty }}.$ Moreover,
 \[{{\left\| {{\left( \lambda -D \right)}^{-1}}\varphi  \right\|}_{1}}\le {{\left\| \varphi  \right\|}_{1}}\int\limits_{0}^{\infty }{{{e}^{-t\left| \operatorname{Re}\lambda  \right|}}dt}=\frac{{{\left\| \varphi  \right\|}_{1}}}{\left| \operatorname{Re}\lambda  \right|}\]
and by complete induction according to $n$ we have
\[{{\left\| {{\left( \lambda -D \right)}^{-n}}\varphi  \right\|}_{1}}\le \frac{{{\left\| \varphi  \right\|}_{1}}}{{{\left| \operatorname{Re}\lambda  \right|}^{n}}}.\]
Excellent method in [4] will show that 
\[\underset{n\to \infty }{\mathop{\lim }}\,\left\| {{\left( \lambda -D \right)}^{-n}}\varphi  \right\|_{1}^{1/n}=\sup \left\{ {{\left| \lambda -i\xi  \right|}^{-1}}:\text{ }\xi \in \text{supp}\left( {\hat{\varphi }} \right) \right\}.\]
It is also proved in [3] that
\[\underset{n\to \infty }{\mathop{\lim }}\,\left\| {(\lambda -D)^{n}}\varphi  \right\|_{1}^{1/n}=\sup \left\{ {{\left| \lambda -i\xi  \right|}}:\text{ }\xi \in \text{supp}\left( {\hat{\varphi }} \right) \right\}.\]
The {\bf Beurling spectrum}  Spec$\left( u \right)$ of  a function $u\in BC\left( \mathbb{R}\to \mathbb{X} \right)$  is defined by
\[\text{Spec}\left( u \right)=\Bigl\{ \xi \in \mathbb{R}:\text{ }\forall \varepsilon >0,\text{ }\exists \varphi \in S:\text{ supp}\hat{\varphi }\subset \left( \xi -\varepsilon ,\xi +\varepsilon  \right),\text{ }\varphi *u\ne 0 \Bigr\}.\]
The Beurling spectral radius $\rho \left( u \right)$ of  $u$  is defined by 
\[\rho \left( u \right)=\sup \Bigl\{ \left| \xi  \right|:\xi \in \text{Spec}\left( u \right) \Bigr\}.\]
For example, let $u\equiv\text{\bf  v}$ (a nonzero vector of $\mathbb{X}$). Then 
\[\varphi *u\left( s \right)=\int\limits_{-\infty }^{\infty }{\varphi \left( s-t \right)dt\text{\bf  v}\equiv \hat{\varphi }\left( 0 \right)}\text{\bf  v}\]
 and consequently, Spec$\left( \text{\bf v} \right)=\left\{ 0 \right\}.$  If $u\equiv 0$ (the zero vector of $\mathbb{X}$) then Spec$\left( 0 \right)=\varnothing .$ Now let $u\left( t \right)=\cos t\text{\bf v }.$ Then 
\[\varphi *u\left( s \right)=\int\limits_{-\infty }^{\infty }{\varphi \left( s-t \right)\cos tdt\text{\bf v }} =\frac{{{e}^{is}}\hat{\varphi }\left( 1 \right)+{{e}^{-is}}\hat{\varphi }\left( -1 \right)}{2}
\cdot\text{\bf  v}\]
and consequently, Spec$\left( \cos t\text{\bf v } \right)=\left\{ 1,-1 \right\}.$ Similarly, 
\[\text{Spec}\left( \text{ }\sum\limits_{k=1}^{n}{{{e}^{it{{\xi }_{k}}}}{\text{\bf v}_{k}}} \right)=\Bigl\{ {{\xi }_{1}},{{\xi }_{2}},\cdots ,{{\xi }_{n}} \Bigr\},\]
where ${{\xi }_{1}},{{\xi }_{2}},\cdots ,{{\xi }_{n}}$ are fixed distinct real numbers and ${\text{\bf v}_{1}},{\text{\bf v}_{2}},\cdots ,{\text{\bf  v}_{n}}$ are fixed nonzero vectors of $\mathbb{X}.$ More generally, if $u\left( t \right)=\phi \left( t \right)\text{\bf  v,}$ where $\phi $ is a bounded  continuous function on the real line, then Spec$\left( u \right)=\text{supp}\left( {\hat{\phi }} \right)$ (the Fourier transform of $\phi $ is taking in the distributional sense). Note that Spec$\left( u \right)$ is always a closed subset of $\mathbb{R}.$
Moreover,
\begin{itemize}
\item $\text{Spec}\left( u+v \right)\subset \text{Spec}\left( u \right)\cup \text{Spec}\left( v \right)$  for all   $u,v\in BC\left( \mathbb{R}\to \mathbb{X} \right);$
\item 	if $u\left( t \right)={{e}^{i\xi t}}v\left( t \right)$ then $\text{Spec}\left( u \right)=\text{Spec}\left( v \right)+\xi ;$
\item 	$\text{Spec}\left( \varphi *u \right)\subset \text{Spec}\left( u \right)\cap \text{supp}\left( {\hat{\varphi }} \right)$ for all $u\in BC\left( \mathbb{R}\to \mathbb{X} \right)$ and $\varphi \in {{L}^{1}}\left( \mathbb{R} \right);$
\item 	if ${{\hat{\varphi }}_{0}}\equiv 0$ on Spec$\left( u \right)$ then ${{\varphi }_{0}}*u=0;$
\item 	if ${{\hat{\varphi }}_{0}}\equiv 1$ on Spec$\left( u \right)$ then ${{\varphi }_{0}}*u=u;$
\item 	$\text{Spec}\left( {{\left( \lambda -D \right)}^{-1}}u \right)=\text{Spec}\left( u \right)$ for every  $u\in BC\left( \mathbb{R}\to \mathbb{X} \right)$ and $\lambda \in \mathbb{C}\backslash i\mathbb{R};$
\item  if ${{\left\{ {{u}_{s}} \right\}}_{s\in \left[ 0,1 \right]}}\subseteq BC\left( \mathbb{R}\to \mathbb{X} \right)$ is a continuous function from $\left[ 0,1 \right]$ into  $BC\left( \mathbb{R}\to \mathbb{X} \right)$ then Spec${{\left( {{u}_{s}} \right)}_{s\in \left[ 0,1 \right]}}$ is a multi-valued continuous function from $\left[ 0,1 \right]$ into ${{2}^{\mathbb{R}}}.$
\end{itemize}
\par\noindent
See \cite{minh} for more details. The following theorem is an analogy of Gelfand famous spectral radius  theorem. It is also an extension of Ha Huy Bang [3]  excellent results for $L^p$ to any Banach space.

\bigskip
\bigskip

\par\noindent {\bf Theorem 1. }  {\it  If  } Spec$\left( u \right)$ {\it  is compact then u is infinitely differentiable, ${{D}^{n}}u\in BC\left( \mathbb{R}\to \mathbb{X} \right)$  for every $n=1,2,\cdots $  and
\[\underset{n\to \infty }{\mathop{\lim }}\,\left\| {{D}^{n}}u \right\|_{\infty }^{1/n}=\rho \left( u \right).\]
Conversely, if u is infinitely differentiable and 
\[\underset{n\to \infty }{\mathop{\lim \inf }}\,\left\| {{D}^{n}}u \right\|_{\infty }^{1/n}<\infty \] then}  Spec$\left( u \right)$ {\it is compact.}

\bigskip
\par\noindent {\sl Proof. } First, assume that $u$  is infinitely differentiable and $\underset{n\to \infty }{\mathop{\lim \inf }}\,\left\| {{D}^{n}}u \right\|_{\infty }^{1/n}<\infty .$  Let ${{\xi }_{0}}\in \text{Spec}\left( u \right)\backslash \left\{ 0 \right\}$ and let $\varepsilon \in \left( 0,\frac{\left| {{\xi }_{0}} \right|}{2} \right).$  According to the definition of Spec$\left( u \right)$ there is a Schwartz function ${{\varphi }_{0}}$ such that $\text{supp}{{\hat{\varphi }}_{0}}\subset \left( {{\xi }_{0}}-\varepsilon ,{{\xi }_{0}}+\varepsilon  \right)$ and ${{\varphi }_{0}}*u\ne 0.$ Let 
\[{{\psi }_{n}}\left( s \right)=\int\limits_{-\infty }^{\infty }{{{e}^{ist }}{{\left( \frac{\left| {{\xi }_{0}} \right|-2\varepsilon }{t} \right)}^{n}}{{{\hat{\varphi }}}_{0}}\left( t \right)dt}.\]
According to [1, p. 507] we have ${{\left\| {{\psi }_{n}} \right\|}_{1}}\le M$ independent of $n$.  Hence,
\[{{\left\| {{D}^{n}}{{\psi }_{n}}*u \right\|}_{\infty }}={{\left\| {{\psi }_{n}}*{{D}^{n}}u \right\|}_{\infty }}\le {{\left\| {{\psi }_{n}} \right\|}_{1}}{{\left\| {{D}^{n}}u \right\|}_{\infty }}\le M{{\left\| {{D}^{n}}u \right\|}_{\infty }}.\]
On the other hand,
\[\begin{aligned}
 {{D}^{n}}{{\psi }_{n}}\left( s \right) 
& =\frac{{{d}^{n}}}{d{{s}^{n}}}\int\limits_{-\infty }^{\infty }{{{e}^{is\xi }}{{\left( \frac{\left| {{\xi }_{0}} \right|-2\varepsilon }{\xi } \right)}^{n}}{{{\hat{\varphi }}}_{0}}\left( \xi  \right)d\xi } \\ 
 & ={{i}^{n}}{{\left( \left| {{\xi }_{0}} \right|-2\varepsilon  \right)}^{n}}\int\limits_{-\infty }^{\infty }{{{e}^{is\xi }}{{{\hat{\varphi }}}_{0}}\left( \xi  \right)d\xi } \\ 
 & ={{i}^{n}}2\pi {{\varphi }_{0}}\left( s \right){{\left( \left| {{\xi }_{0}} \right|-2\varepsilon  \right)}^{n}},  
\end{aligned}\]
and consequently,
\[2\pi {{\left\| {{\varphi }_{0}}*u \right\|}_{\infty }}{{\left( \left| {{\xi }_{0}} \right|-2\varepsilon  \right)}^{n}}\le M{{\left\| {{D}^{n}}u \right\|}_{\infty }}.\]
Since ${{\varphi }_{0}}*u\ne 0,$
\[\left| {{\xi }_{0}} \right|-2\varepsilon \le \underset{n\to \infty }{\mathop{\lim \inf }}\,\left\| {{D}^{n}}u \right\|_{\infty }^{1/n}.\]
But ${{\xi }_{0}}\in \text{Spec}\left( u \right)\backslash \left\{ 0 \right\}$ is arbitrary, so 
\[\rho \left( u \right)\le \underset{n\to \infty }{\mathop{\lim \inf }}\,\left\| {{D}^{n}}u \right\|_{\infty }^{1/n}<\infty \]
and Spec$\left( u \right)$ is compact. Conversely, assume that Spec$\left( u \right)$ is compact. Let ${{\varphi }_{0}}$ be a Schwartz function such that supp$\left( {{{\hat{\varphi }}}_{0}} \right)\subset \left[ -\rho \left( u \right)-\varepsilon ,\rho \left( u \right)+\varepsilon  \right]$ and ${{\hat{\varphi }}_{0}}\equiv 1$ on Spec$\left( u \right)$ ($\varepsilon >0$ given). Then  $u={{\varphi }_{0}}*u,$ and consequently, ${{D}^{n}}u={{D}^{n}}{{\varphi }_{0}}*u\in BC\left( \mathbb{R}\to \mathbb{X} \right)$ for every $n=1,2,\cdots .$ Hence, $u$  is infinitely differentiable and $\rho \left( u \right)\le \underset{n\to \infty }{\mathop{\lim \inf }}\,\left\| {{D}^{n}}u \right\|_{\infty }^{1/n}.$ Moreover,
\[\begin{aligned}
\underset{n\to \infty }{\mathop{\lim \sup }}\,\left\| {{D}^{n}}u \right\|_{\infty }^{1/n}
&=\underset{n\to \infty }{\mathop{\lim \sup }}\,\left\| {{D}^{n}}{{\varphi }_{0}}*u \right\|_{\infty }^{1/n}\\
&\le \underset{n\to \infty }{\mathop{\lim \sup }}\,\left\| {{D}^{n}}{{\varphi }_{0}} \right\|_{1}^{1/n}\left\| u \right\|_{\infty }^{1/n}\\
&\le \rho \left( u \right)+\varepsilon
\end{aligned}\]
(by  [1, p. 506]). 
Thus, $\underset{n\to \infty }{\mathop{\lim \sup }}\,\left\| {{D}^{n}}u \right\|_{\infty }^{1/n}\le \rho \left( u \right).$ But we have proved that $\rho \left( u \right)\le \underset{n\to \infty }{\mathop{\lim \inf }}\,\left\| {{D}^{n}}u \right\|_{\infty }^{1/n}$
 so the proof is now complete.

\bigskip

\par\noindent{\bf  Remark. } If $K\subseteq \mathbb{R}$ is a compact set and \[\mathbb{V}\left( K \right)=\left\{ u\in BC\left( \mathbb{R}\to \mathbb{X} \right):\text{ Spec}\left( u \right)\subseteq K \right\}\]
 then the differential operator D  is bounded in $\mathbb{V}\left( K \right)$ and Spec$\left( {{\left. D \right|}_{\mathbb{V}\left( K \right)}} \right)=iK$ and the usual spectral theory can be applied.

\bigskip

\par\noindent {\bf Corollary. }  {\it If }  Spec$\left( \text{ }u \right)=\left\{ {{\xi }_{1}},{{\xi }_{2}},\cdots ,{{\xi }_{n}} \right\},$ {\it then $u\left( t \right)=\sum\limits_{k=1}^{n}{{{e}^{it{{\xi }_{k}}}}{{\mathbf v}_{k}}}$ where ${{\mathbf v}_{1}},{{\mathbf v}_{2}},\cdots ,{{\mathbf  v}_{n}}$ are fixed nonzero vectors of $\mathbb{X}.$}

\bigskip

\par\noindent{\sl Proof. }  For a real number $\xi $ let 
\[\mathbb{V}\left( {{\xi }_{1}},{{\xi }_{2}},\cdots ,{{\xi }_{n}} \right)=\left\{ u\in BC\left( \mathbb{R}\to \mathbb{X} \right):\text{ Spec}\left( u \right)\subseteq \left\{ {{\xi }_{1}},{{\xi }_{2}},\cdots ,{{\xi }_{n}} \right\} \right\}.\]
By the above theorem,
\[\begin{aligned}
  \mathbb{V}\left( 0 \right) &=\left\{ u\in BC\left( \mathbb{R}\to \mathbb{X} \right):\text{ Spec}\left( u \right)\subseteq \left\{ 0 \right\} \right\}\\
&=\left\{ u\in BC\left( \mathbb{R}\to \mathbb{X} \right):\text{ }Du=0 \right\} \\ 
 & =\left\{ \text{constant functions} \right\}  
\end{aligned}\]
hence, 
\[\mathbb{V}\left( \xi  \right)=\Bigl\{ u\in BC\left( \mathbb{R}\to \mathbb{X} \right):\text{ }u\left( t \right)={{e}^{i\xi t}}{\mathbf v}\quad \text{\bf v}\in \mathbb{X} \Bigr\}.\]
Note that 
\[\mathbb{V}\left( {{\xi }_{1}},{{\xi }_{2}},\cdots ,{{\xi }_{n}} \right)=\underset{j=1}{\overset{n}{\mathop{\oplus }}}\,\mathbb{V}\left( {{\xi }_{j}} \right)\]
and this complete the proof. (Note that this corollary was proved by N.V. Minh \cite{minh} by other way).

\bigskip

\par\noindent{\bf Example 1.}  Let $A$ be a bounded linear operator on $\mathbb{X}$ and put $u\left( t \right)={{e}^{iAt}}{\mathbf v}$ ({\bf v } is a nonzero vector of $\mathbb{X}$). Then ${{D}^{n}}u={{\left( iA \right)}^{n}}u,$ so  $\rho \left( u \right)\le \rho \left( A \right).$

\bigskip

\par\noindent{\bf Example 2.}  Consider the delay equation $\dot{u}\left( t \right)=-u\left( t-\tau  \right)$  for all  $t\in \mathbb{R}.$  We can easily compute ${{D}^{n}}u\left( t \right)={{\left( -1 \right)}^{n}}u\left( t-n\tau  \right)$  and get  $\underset{n\to \infty }{\mathop{\lim }}\,\left\| {{D}^{n}}u \right\|_{\infty }^{1/n}=1$  if $u$  is bounded and non-identically 0.  Hence, in this case  $\rho \left( u \right)=1.$

\bigskip

The following theorem is an extension of [2] and [4]. In proof we will use the method of Ha Huy Bang  [4].

\bigskip
\par\noindent{\bf Theorem 2.}  {\it   If  $u\in BC\left( \mathbb{R}\to \mathbb{X} \right)$ and $\lambda \in \mathbb{C}\backslash i\mathbb{R}$ then}
\[\underset{n\to \infty }{\mathop{\lim \inf }}\,\left\| {{\left( \lambda -D \right)}^{-n}}u \right\|_{\infty }^{1/n}\ge \sup \left\{ {{\left| \lambda -i\xi  \right|}^{-1}}:\text{ }\xi \in \text{Spec}\left( u \right) \right\}.\]
{\it Moreover, if  }  Spec$\left( u \right)$  {\it is compact then}
\[\underset{n\to \infty }{\mathop{\lim }}\,\left\| {{\left( \lambda -D \right)}^{-n}}u \right\|_{\infty }^{1/n}=\sup \left\{ {{\left| \lambda -i\xi  \right|}^{-1}}:\text{ }\xi \in \text{Spec}\left( u \right) \right\}.\]

\bigskip
\par\noindent{Proof. } Let $\xi \in \text{Spec}\left( u \right).$ According to the definition of Spec$\left( u \right),$ for every $\varepsilon >0,$ there is a Schwartz function $\varphi $ such that $\text{supp}\hat{\varphi }\subset \left( \xi -\varepsilon ,\xi +\varepsilon  \right)$ and $\varphi *u\ne 0.$ Note that 
\[\varphi \left( s \right)=\frac{1}{2\pi }\int\limits_{\mathbb{R}}{{{e}^{ist}}\overset{\lower0.5em\hbox{$\smash{\scriptscriptstyle\frown}$}}{\varphi}\left(t\right)dt}\quad\text{ and }\quad {{\left( \lambda -D \right)}^{n}}\varphi \left( s \right)=\frac{1}{2\pi }\int\limits_{\mathbb{R}}{{{e}^{ist}}{{\left( \lambda -it \right)}^{n}}\overset{\lower0.5em\hbox{$\smash{\scriptscriptstyle\frown}$}}{\varphi}\left(t\right)dt}.\]
Let $a=\underset{t\in \left( \xi -\varepsilon ,\xi +\varepsilon  \right)}{\mathop{\text{sup}}}\,\left| \lambda -it \right|+\varepsilon $ and 
\[{{\psi }_{n}}\left( s \right)={{a}^{-n}}\int\limits_{\mathbb{R}}{{{e}^{ist}}{{\left( \lambda -it \right)}^{n}}\overset{\lower0.5em\hbox{$\smash{\scriptscriptstyle\frown}$}}{\varphi}\left(t\right)dt}=2\pi{{a}^{-n}}{{\left(\lambda-D\right)}^{n}}\varphi\left(s\right).\]
We prove that  ${{\left\| {{\psi }_{n}} \right\|}_{1}}\le M$ independent of $n$.  
Indeed, 
\[\begin{aligned}
  {{\left\| {{\psi }_{n}} \right\|}_{1}}
 &\le \int\limits_{-\infty }^{\infty }{\frac{dx}{{{x}^{2}}+1}}\cdot \underset{x\in \mathbb{R}}{\mathop{\sup }}\,\left| \left( {{x}^{2}}+1 \right){{\psi }_{n}}\left( x \right) \right| \\ 
 & =2{{\pi }^{2}}{{a}^{-n}}\underset{x\in \mathbb{R}}{\mathop{\sup }}\,\left| \left( {{x}^{2}}+1 \right){{\left( \lambda -D \right)}^{n}}\varphi \left( x \right) \right| \\ 
  & \le Ca^{-n}n^2(a-\epsilon )^n\\
 & \le M\text{ independent of }n.  
\end{aligned}\]
Hence,
\[\begin{aligned}{{\left\| {{\left( \lambda -D \right)}^{-n}}{{\psi }_{n}}*u \right\|}_{\infty }}
&={{\left\| {{\psi }_{n}}*{{\left( \lambda -D \right)}^{-n}}u \right\|}_{\infty }}\\
&\le {{\left\| {{\psi }_{n}} \right\|}_{1}}{{\left\| {{\left( \lambda -D \right)}^{-n}}u \right\|}_{\infty }}\\
&\le M{{\left\| {{\left( \lambda -D \right)}^{-n}}u \right\|}_{\infty }}.
\end{aligned}\]
On the other hand,
\[{{\left( \lambda -D \right)}^{-n}}{{\psi }_{n}}\left( s \right)=2\pi {{a}^{-n}}{{\left( \lambda -D \right)}^{n}}\varphi \left( s \right)=2\pi \varphi \left( s \right){{a}^{-n}},\]
and consequently,
\[2\pi {{\left\| \varphi *u \right\|}_{\infty }}{{a}^{-n}}\le M{{\left\| {{\left( \lambda -D \right)}^{-n}}u \right\|}_{\infty }}.\]
Since $\varphi *u\ne 0,$
\[{{a}^{-1}}\le \underset{n\to \infty }{\mathop{\lim \inf }}\,\left\| {{\left( \lambda -D \right)}^{-n}}u \right\|_{\infty }^{1/n}.\]
But $\xi \in \text{Spec}\left( u \right)$ and $\varepsilon >0$ are arbitrary, so 
\[\sup \left\{ {{\left| \lambda -i\xi  \right|}^{-1}}:\text{ }\xi \in \text{Spec}\left( u \right) \right\}\le \underset{n\to \infty }{\mathop{\lim \inf }}\,\left\| {{\left( \lambda -D \right)}^{-n}}u \right\|_{\infty }^{1/n}.\]
Next we prove that
\[\underset{n\to \infty }{\mathop{\lim \sup }}\,\left\| {{\left( \lambda -D \right)}^{-n}}u \right\|_{\infty }^{1/n}\le \sup \left\{ {{\left| \lambda -i\xi  \right|}^{-1}}:\text{ }\xi \in \text{Spec}\left( u \right) \right\}.\]
if  Spec$\left( u \right)$ is compact. Indeed, if Spec$\left( u \right)$ is compact then $u=\varphi *u$ for a Schwartz function $\varphi $ ($\overset{\lower0.5em\hbox{$\smash{\scriptscriptstyle\frown}$}}{\varphi}=1$ on the spectrum of $u$). Consequently,  
\[{{\left( \lambda -D \right)}^{-n}}u={{\left( \lambda -D \right)}^{-n}}\varphi *u
\]
 and  
\[\underset{n\to \infty }{\mathop{\lim \sup }}\,\left\| {{\left( \lambda -D \right)}^{-n}}u \right\|_{\infty }^{1/n}\le \underset{n\to \infty }{\mathop{\lim \sup }}\,\left\| {{\left( \lambda -D \right)}^{-n}}\varphi  \right\|_{1}^{1/n}.\] On the other hand, 
\[{{\left( \lambda -D \right)}^{-n}}\varphi \left( s \right)=\frac{1}{2\pi }\int\limits_{\mathbb{R}}{{{e}^{ist}}{{\left( \lambda -it \right)}^{-n}}\overset{\lower0.5em\hbox{$\smash{\scriptscriptstyle\frown}$}}{\varphi}\left(t\right)dt}\]
and similar to [4] we have
\[\begin{aligned}
 {{\left\| {{\left( \lambda -D \right)}^{-n}}\varphi  \right\|}_{1}}
 & =\underset{{{\left\| \psi  \right\|}_{\infty }}\le 1}{\mathop{\sup }}\,\left| \int\limits_{-\infty }^{\infty }{{{\left( \lambda -D \right)}^{-n}}\varphi \left( s \right)\psi \left( s \right)ds} \right| \\ 
 & =\frac{1}{2\pi }\underset{{{\left\| \psi  \right\|}_{\infty }}\le 1}{\mathop{\sup }}\,\left| \int\limits_{-\infty }^{\infty }{\overset{\lower0.5em\hbox{$\smash{\scriptscriptstyle\frown}$}}{\varphi}\left(t\right){{\left(\lambda-it\right)}^{-n}}dt}\int\limits_{-\infty}^{\infty}{{{e}^{ist}}\psi\left(s\right)ds}\right|\\&=\underset{{{\left\|\psi\right\|}_{\infty}}\le1}{\mathop{\sup}}\,\left|\int\limits_{-\infty}^{\infty}{\hat{\varphi}\left(t\right){{\left(\lambda-it\right)}^{-n}}h\left(t\right){{\cal F}^{-1}}\left(\psi,t\right)dt}\right|\\
&=\underset{{{\left\|\psi\right\|}_{\infty}}\le1}{\mathop{\sup}}\,\left|\int\limits_{-\infty}^{\infty}{\varphi\left(t\right){\cal F}\left(\frac{h{{\cal F}^{-1}}\psi}{{{\left(\lambda-i\tau\right)}^{n}}},t\right)dt}\right|\\&=\frac{1}{2\pi}\underset{{{\left\|\psi\right\|}_{\infty}}\le1}{\mathop{\sup}}\,\left|\int\limits_{-\infty}^{\infty}{\varphi\left(t\right)\left({\cal F}\left(\frac{h}{{{\left(\lambda-i\tau\right)}^{n}}}\right)*{\cal F}\left({{\cal F}^{-1}}\psi\right)\right)\left(t\right)dt}\right|\\
&=\frac{1}{2\pi}\underset{{{\left\|\psi\right\|}_{\infty}}\le1}{\mathop{\sup}}\,\left|\int\limits_{-\infty}^{\infty}{\varphi\left(t\right)\left({\cal F}\left(\frac{h}{{{\left(\lambda-i\tau\right)}^{n}}}\right)*\psi\right)\left(t\right)dt}\right|\\&\le\frac{1}{2\pi}\underset{{{\left\|\psi\right\|}_{\infty}}\le1}{\mathop{\sup}}\,{{\left\|\varphi\right\|}_{1}}{{\left\|{\cal F}\left(\frac{h}{{{\left(\lambda-i\tau\right)}^{n}}}\right)*\psi\right\|}_{\infty}}\\
&\le \frac{{{\left\| \varphi  \right\|}_{1}}}{2\pi }{{\left\| {\cal F}\left( \frac{h}{{{\left( \lambda -i\tau  \right)}^{n}}} \right) \right\|}_{1}},
\end{aligned}\]
where $h$  is a test function such that $h=1$ on the support of $\overset{\lower0.5em\hbox{$\smash{\scriptscriptstyle\frown}$}}{\varphi}.$ We assume that $h\left( \tau  \right)=0$ for $\left| \tau  \right|>\rho \left( u \right)+\varepsilon .$ We have 
\[\begin{aligned}
   {\cal F}\left( \frac{h}{{{\left( \lambda -it \right)}^{n}}},s \right)
&=\int\limits_{\left| \tau  \right|<\rho \left( u \right)+\varepsilon .}{\frac{{{e}^{-its}}h\left( t \right)}{{{\left( \lambda -it \right)}^{n}}}dt}=\int\limits_{\left| \tau  \right|<\rho \left( u \right)+\varepsilon .}{\frac{{{\left( {{e}^{-its}} \right)}^{\prime }}h\left( t \right)}{-is{{\left( \lambda -it \right)}^{n}}}dt} \\ 
 & =\int\limits_{\left| \tau  \right|<\rho \left( u \right)+\varepsilon .}{\frac{{{e}^{-its}}{{h}^{\prime }}\left( t \right)}{is{{\left( \lambda -it \right)}^{n}}}dt}+\int\limits_{\left| \tau  \right|<\rho \left( u \right)+\varepsilon .}{\frac{n{{e}^{-its}}h\left( t \right)}{s{{\left( \lambda -it \right)}^{n+1}}}dt}.  
\end{aligned}\]
Therefore,  
\[\begin{aligned}
 \underset{n\to \infty }{\mathop{\lim \sup }}\,\left\| {{\left( \lambda -D \right)}^{-n}}\varphi  \right\|_{1}^{1/n} 
& \le \sup \left\{ {{\left| \lambda -i\xi  \right|}^{-1}}:\text{ }\xi \in \text{supp}\left( {\overset{\lower0.5em\hbox{$\smash{\scriptscriptstyle\frown}$}}{\varphi}}\right)\right\}\\
&\le\varepsilon+\sup\left\{{{\left|\lambda-i\xi\right|}^{-1}}:\text{}\xi\in\text{Spec}\left(u\right)\right\}\end{aligned}\]
and  this completes the proof.

\bigskip
\bigskip

\sect{\bf  Fourier coefficients of almost periodic functions
 and the Beurling spectrum}

\bigskip

\par\noindent Now we are interested in computing the spectrum of almost periodic functions. To this end we define the $\lambda $th Fourier coefficient ${{a}_{\lambda }}\left( u \right)$ (a vector  of $\mathbb{X}$)  of  a function $u\in BC\left( \mathbb{R}\to \mathbb{X} \right)$ by letting
\[{{a}_{\lambda }}\left( u \right)=\underset{T\to \infty }{\mathop{\lim }}\,\frac{1}{2T}\int\limits_{-T}^{T}{{{e}^{-i\lambda t}}u\left( t \right)dt.}\]
If this limit exists for every $\lambda \in \mathbb{R},$ the function $u\in BC\left( \mathbb{R}\to \mathbb{X} \right)$ is called almost periodic.  It is easy to prove that for every almost periodic function $u\in BC\left( \mathbb{R}\to \mathbb{X} \right),$ $\text{Spec}\left( u \right)$ is the closure of $\left\{ \lambda \in \mathbb{R}:\text{ }{{a}_{\lambda }}\left( u \right)\ne 0 \right\}.$ For  example, let $u\left( t \right)={{e}^{it}}\text{\bf v }.$ We have
\[\begin{aligned}
{{a}_{\lambda }}\left( u \right)  & =\underset{T\to \infty }{\mathop{\lim }}\,\frac{1}{2T}\int\limits_{-T}^{T}{{{e}^{-i\left( \lambda -1 \right)t}}dt}\cdot {\mathbf  v}= {\mathbf  v}\quad\text{ if } \lambda=1, \\ 
 & =\underset{T\to \infty }{\mathop{\lim }}\,\frac{\sin \left( \lambda -\text{1} \right)T}{\left( \lambda -\text{1} \right)T}\cdot {\mathbf v}=0\quad\text{  if  }\lambda \ne 1,  
\end{aligned}\]
and consequently, $\text{Spec}\left( {{e}^{it}}\text{\bf v }\right)=\left\{ 1 \right\}.$ More generally, let $u\left( t \right)={{e}^{iAt}}{\mathbf v}$ ({\bf v}  is a nonzero vector of $\mathbb{X}$ and  $A:\mathbb{X}\to \mathbb{X}$ is a bounded linear operator on $\mathbb{X}$). Then we have $\dot{u}\left( t \right)=iAu\left( t \right),$ so if $u$ is almost periodic, 
\[iA{{a}_{\lambda }}\left( u \right)=\underset{T\to \infty }{\mathop{\lim }}\,\frac{1}{2T}\int\limits_{-T}^{T}{{{e}^{-i\lambda t}}\dot{u}\left( t \right)dt}=\underset{T\to \infty }{\mathop{\lim }}\,\frac{1}{2T}\int\limits_{-T}^{T}{i\lambda {{e}^{-i\lambda t}}u\left( t \right)dt}=i\lambda {{a}_{\lambda }}\left( u \right).
\]
Hence, if $\lambda \in \text{Spec}\left( u \right)$ then $\lambda $ is a real eigenvalue of $A$, and consequently, Spec$\left( u \right)\subseteq \text{Spec}\left( A \right)$ (the point spectrum of $A$). Now consider the delay differential equation [5] [7]
$\dot{u}\left( t \right)=-u\left( t-\tau  \right)$  for all $t\in {\mathbb R}.$ If $u$ is almost periodic, the Fourier coefficient ${{a}_{\lambda }}\left( u \right)$ of $u$ is satisfying $\left( 1+i\lambda {{e}^{i\lambda \tau }} \right){{a}_{\lambda }}\left( u \right)=0.$ Hence, $\text{Spec}\left( u \right)\subseteq \left\{ \lambda \in \left[ -1,1 \right]:\text{ }1+i\lambda {{e}^{i\lambda \tau }}=0 \right\}.$  Therefore, if $u\ne 0,$ we conclude that the equation $1+i\lambda {{e}^{i\lambda \tau }}=0$ should have a real root in $\left[ -1,1 \right].$ This implies that $\tau =\frac{\pi }{2}$ and Spec$\left( u \right)\subseteq \left\{ \pm 1 \right\}.$ 
(See [5] for more details.)
 Hence, $u\left( t \right)=$ $\cos t\text{ }{{\mathbf v}_{1}}+\sin t\text{ }{{\mathbf v}_{2}}$  is periodic (${{\mathbf v}_{1}}\text{ and  }{{\mathbf v}_{2}}$ are vectors in the Banach space $\mathbb{X}$). Hence the delay equation $\dot{u}\left( t \right)=-u\left( t-\tau  \right)$  has a non-zero almost periodic solution if and only if $\tau =\frac{\pi }{2}$. Now consider the delay equation $\dot{u}\left( t \right)=iAu\left( t-\tau  \right)$ for all $t\in \mathbb{R},$  ($A:\mathbb{X}\to \mathbb{X}$ is bounded and linear).  We can compute easily $\rho \left( u \right)\le \rho \left( A \right)$ (the spectral radius of $A$). Moreover, if $u$ is almost periodic, the Fourier coefficient ${{a}_{\lambda }}\left( u \right)$ satisfies $A{{a}_{\lambda }}\left( u \right)=\lambda {{e}^{i\lambda \tau }}{{a}_{\lambda }}\left( u \right).$ Hence, if $\lambda \in \text{Spec}\left( u \right),$ then $\lambda {{e}^{i\lambda \tau }}$ is an eigenvalue of $A$. Therefore, $\text{Spec}\left( u \right)\subseteq $ $\left\{ \lambda \in \left[ -\rho \left( A \right),\rho \left( A \right) \right]:\text{ }\lambda {{e}^{i\lambda \tau }}\in \text{Spec}\left( A \right) \right\}.$  Specially, if $\tau =0$ we get back  the above result Spec$\left( u \right)\subseteq \text{Spec}\left( A \right)\cap \mathbb{R}.$

\bigskip

\par\noindent{\bf Acknowledgement.} The author would like to express his sincerely thank to the referee for reading carefully the manuscript and providing some suggestions  that have been implemented in the final version of the paper. 
Deepest appreciation is extended towards the NAFOSTED  (the National Foundation for Science and Techology Development in Vietnam) for the financial support.

\end{document}